%% file: sketch_alessandra.tex
\documentclass[twoside]{aiml16}
\usepackage{aiml16macro}
\usepackage{cmll}
\usepackage{graphicx,verbatim, bussproofs,tikz}
\usepackage{pgfplots}
\usepackage{amsmath}
\usepackage{amssymb}
\usepackage{txfonts}
\usepackage{graphicx}
\usetikzlibrary{matrix,arrows}
\usetikzlibrary{patterns}
\usetikzlibrary{shapes}

\tikzset{
	treenode/.style = {align=center, inner sep=0pt, text centered},
	Ske/.style = {treenode, ellipse, double, draw=black,
		minimum width=6pt, thick},
	PIA/.style = {treenode, ellipse, black, draw=black,
		minimum width=6pt},
	Crit/.style = {treenode, rectangle, draw=black,
		minimum width=0.5em, minimum height=0.5em}
}
\renewcommand{\phi}{\varphi}

\font\bmi=cmmi8 scaled 1440
\newcommand{\powerset}{\raise.6ex\hbox{\bmi\char'175 }}
\font\intix=cmex9

\newcommand{\TM}{{{\hbox{$\textfont3=\intix{\bigwedge}$}}}}
\newcommand{\TJ}{{{{\hbox{$\textfont3=\intix{\bigvee}$}}}}}
\newcommand{\NOM}{\mathsf{NOM}}
\newcommand{\CNOM}{\mathsf{CNOM}}

\newcommand{\ox}{\overline{x}}
\newcommand{\oy}{\overline{y}}
\newcommand{\oz}{\overline{z}}
\newcommand{\oX}{\overline{X}}

\newcommand{\ophi}{\overline{\phi}}
\newcommand{\opsi}{\overline{\psi}}
\newcommand{\oga}{\overline{\gamma}}

\newcommand{\ifBox}[1]{IF$^{\Box}_{#1}$}
\newcommand{\ifDia}[1]{IF$^{\Diamond}_{#1}$}

\newcommand{\LFPtwo}{\mu_{2}}
\newcommand{\LFPstar}{\mu^{*}}

\newcommand{\GFPtwo}{\nu_{2}}
\newcommand{\Tm}{\mathcal{L}}
\newcommand{\Tmone}{\mathcal{L}_1}
\newcommand{\Tmtwo}{\mathcal{L}_2}
\newif\ifmargincoms

\newcommand{\commment}[1]{}

\renewcommand{\phi}{\varphi}
\renewcommand{\emptyset}{\varnothing}
\newcommand{\ca}{\mathbb{A}^\delta}
\newcommand{\A}{\mathbb{A}}

\newcommand{\jira}{J^{\infty}(\bba)}
\newcommand{\mira}{M^{\infty}(\bba)}

\renewcommand{\epsilon}{\varepsilon}

\newcommand{\nomi}{\mathbf{i}}
\newcommand{\nomj}{\mathbf{j}}

\newcommand{\cnomm}{\mathbf{m}}
\newcommand{\cnomn}{\mathbf{n}}

\newcommand{\bba}{\mathbb{A}}

\newcommand{\bbA}{\mathbb{A}}

\newcommand{\bbas}{\bbA^{\delta}}

\newcommand{\kbbas}{K(\mathbb{A}^\delta)}
\newcommand{\obbas}{O(\mathbb{A}^\delta)}

\newcommand{\jty}{J^{\infty}}
\newcommand{\mty}{M^{\infty}}

\newcommand{\bigamp}{\mathop{\mbox{\Large \&}}}

\newcommand{\amp}{\mathop{\&}}
\newcommand{\marginnote}[1]{\marginpar{\raggedright\tiny{#1}}}

\def\lastname{Conradie, Craig, Palmigiano and Zhao}

\begin{document}

\begin{frontmatter} \title{Constructive canonicity for lattice-based fixed point logics}
 \author{Willem Conradie}\footnote{The research of the first author has been made possible by the National Research Foundation of South Africa, Grant number 81309.}
 \author{Andrew Craig}
 \address{Department of Pure and Applied Mathematics, University of Johannesburg, South Africa}
 \author{Alessandra Palmigiano}\footnote{The research of the third and fourth author has been made possible by the NWO Vidi grant 016.138.314, by the NWO Aspasia grant 015.008.054, and by a Delft Technology Fellowship awarded in 2013.}
 \address{Department of Pure and Applied Mathematics, University of Johannesburg, South Africa}
 \address{Faculty of Technology, Policy and Management, Delft University of Technology, the Netherlands}
 \author{Zhiguang Zhao}
 \address{Faculty of Technology, Policy and Management, Delft University of Technology, the Netherlands}

\begin{abstract}
We prove the algorithmic canonicity of two classes of $\mu$-inequalities in a constructive meta-theory of normal lattice expansions. This result simultaneously generalizes Conradie and Craig's canonicity for $\mu$-inequalities based on a bi-intuitionistic bi-modal language, \cite{CoCr14},
and Conradie and Palmigiano's constructive canonicity for inductive inequalities \cite{ConPal13}
(restricted to normal lattice expansions to keep the page limit). Besides the greater generality, the unification of these strands smoothes the existing treatments for the canonicity of $\mu$-formulas and inequalities. In particular, the rules of the algorithm ALBA used for this result have exactly the same formulation as those of \cite{ConPal13}, with no additional rule added specifically to handle the fixed point binders. Rather, fixed points are accounted for by certain restrictions on the application of the rules, concerning the order-theoretic properties of the term functions associated with the formulas to which the rules are applied.   
\end{abstract}

\end{frontmatter}
\section{Introduction}

\input{introduction}

\section{Preliminaries}\label{Sec:Prelim}

In this section we collect the needed preliminaries about the syntax and algebraic semantics of lattice-based mu-calculi.  

\subsection{The language of lattice expansions}

In the present subsection, we introduce the propositional fragments of the lattice-based $\mu$-languages we consider in this paper. 
We will make use of the following auxiliary definitions: an {\em order-type} over $n\in \mathbb{N}$\footnote{Throughout the paper, order-types will be typically associated with arrays of variables $\overline{p}: = (p_1,\ldots, p_n)$. 
} is an element $\epsilon\in \{1, \partial\}^n$. For every order-type $\epsilon$, we denote its {\em opposite} order-type by $\epsilon^\partial$, that is, $\epsilon^\partial_i = 1$ iff $\epsilon_i=\partial$ for every $1 \leq i \leq n$. For any lattice $\bba$, we let $\bba^1: = \bba$ and $\bba^\partial$ be the dual lattice, that is, the lattice associated with the converse partial order of $\bba$. For any order-type $\varepsilon$, we let $\bba^\varepsilon: = \Pi_{i = 1}^n \bba^{\varepsilon_i}$. Sometimes we will  write $\top^1$ and $\top^\partial$ for $\top$ and $\bot$ respectively. Similarly, we will write $\bot^1$ and $\bot^\partial$ for $\bot$ and $\top$ respectively. For both order-types and tuples of variables, we will use the symbol $\oplus$ to denote concatenation.

The LE-language $\mathcal{L}(\mathcal{F}, \mathcal{G})$ (abbreviated as $\mathcal{L}$ when no confusion arises) takes as parameters: 1) a denumerable set $\mathsf{PROP}$ of proposition letters, elements of which are denoted $p,q,r$, possibly with indices; 2) disjoint sets of connectives $\mathcal{F}$ and $\mathcal{G}$. Each $f\in \mathcal{F}$ (resp.\ $g\in \mathcal{G}$) has arity $n_f\in \mathbb{N}$ (resp.\ $n_g\in \mathbb{N}$), and is associated with some order-type $\varepsilon_f$ over $n_f$ (resp.\ $\varepsilon_g$ over $n_g$).
The terms (formulas) of $\mathcal{L}$ are defined recursively as follows:
\[\phi ::= p \mid \bot \mid \top \mid \phi \wedge \phi \mid \phi \vee \phi \mid f(\overline{\phi}) \mid g(\overline{\phi})\]
where $p \in \mathsf{PROP}$, $f \in \mathcal{F}$, $g \in \mathcal{G}$. Terms in $\mathcal{L}$ will be denoted either by $s,t$, or by lowercase Greek letters such as $\varphi, \psi, \gamma$ etc.

\subsection{Lattice expansions, and their canonical extensions}
The following definition captures the algebraic setting of the present paper, namely the normal lattice expansions of \cite{ConPal13}.

\begin{definition}
\label{def:DLE}
For any LE-signature $\mathcal{L} = \mathcal{L}(\mathcal{F}, \mathcal{G})$, a {\em lattice expansion} (LE) is a tuple $\bba = (L, \mathcal{F}^\bbA, \mathcal{G}^\bbA)$ such that $L$ is a bounded lattice, $\mathcal{F}^\bbA = \{f^\bbA\mid f\in \mathcal{F}\}$ and $\mathcal{G}^\bbA = \{g^\bbA\mid g\in \mathcal{G}\}$, such that every $f^\bbA\in\mathcal{F}^\bbA$ (resp.\ $g^\bbA\in\mathcal{G}^\bbA$) is an $n_f$-ary (resp.\ $n_g$-ary) operation on $\bbA$. An LE $\bba$ is {\em normal} if every $f^\bbA\in\mathcal{F}^\bbA$ (resp.\ $g^\bbA\in\mathcal{G}^\bbA$) preserves finite (hence also empty) joins (resp.\ meets) in each coordinate with $\epsilon_f(i)=1$ (resp.\ $\epsilon_g(i)=1$) and reverses finite (hence also empty) meets (resp.\ joins) in each coordinate with $\epsilon_f(i)=\partial$ (resp.\ $\epsilon_g(i)=\partial$). Let $\mathbb{LE}$ be the class of LEs.
Sometimes we will refer to certain LEs as $\mathcal{L}$-algebras to emphasize the signature. In the remainder of the paper, we will abuse notation and write e.g.\ $f$ for $f^\bbA$ when this causes no confusion. From now on, normal lattice expansions  will be abbreviated as LEs.
\end{definition}
Each language $\mathcal{L}$ is interpreted in the appropriate class of LEs. In particular, for every LE $\bba$, each operation $f^\bba\in \mathcal{F}^\bbA$ (resp.\ $g^\bba\in \mathcal{G}^\bbA$) is finitely join-preserving (resp.\ meet-preserving) in each coordinate when regarded as a map $f^\bba: \bba^{\varepsilon_f}\to \bba$ (resp.\ $g^\bba: \bba^{\varepsilon_g}\to \bba$).


\subsection{The `tense' language $\mathcal{L}^t$}\label{ssec:expanded tense language}
Any language $\mathcal{L} = \mathcal{L}(\mathcal{F}, \mathcal{G})$ can be associated with the language $\mathcal{L}^t = \mathcal{L}_\mathrm{LE}(\mathcal{F}^t, \mathcal{G}^t)$, where $\mathcal{F}^t\supseteq \mathcal{F}$ and $\mathcal{G}^t\supseteq \mathcal{G}$ are obtained by adding:
\begin{enumerate}
\item for $f\in\mathcal{F}$ and $1\leq i\leq n_f$, the $n_f$-ary connective $f^\sharp_i$, the intended interpretation of which is the right residual of $f$ in its $i$th coordinate if $\varepsilon_f(i) = 1$ (resp.\ its Galois-adjoint if $\varepsilon_f(i) = \partial$);
\item for $g\in\mathcal{G}$ and $1 \leq i\leq n_g$, the $n_g$-ary connective $g^\flat_i$, the intended interpretation of which is the left residual of $g$ in its $i$th coordinate if $\varepsilon_g(i) = 1$ (resp.\ its Galois-adjoint if $\varepsilon_g(i) = \partial$).
\end{enumerate}
We stipulate that $f^\sharp_i\in\mathcal{G}^t$ if $\varepsilon_f(i) = 1$, and $f^\sharp_i\in\mathcal{F}^t$ if $\varepsilon_f(i) = \partial$. Dually, $g^\flat_i\in\mathcal{F}^t$ if $\varepsilon_g(i) = 1$, and $g^\flat_i\in\mathcal{G}^t$ if $\varepsilon_g(i) = \partial$. The order-type assigned to the additional connectives is predicated on the order-type of their intended interpretations. That is, for any $f\in \mathcal{F}$ and $g\in\mathcal{G}$,
\begin{enumerate}
\item if $\epsilon_f(i) = 1$, then $\epsilon_{f_i^\sharp}(i) = 1$ and $\epsilon_{f_i^\sharp}(j) = (\epsilon_f(j))^\partial$ for any $j\neq i$.
\item if $\epsilon_f(i) = \partial$, then $\epsilon_{f_i^\sharp}(i) = \partial$ and $\epsilon_{f_i^\sharp}(j) = \epsilon_f(j)$ for any $j\neq i$.
\item if $\epsilon_g(i) = 1$, then $\epsilon_{g_i^\flat}(i) = 1$ and $\epsilon_{g_i^\flat}(j) = (\epsilon_g(j))^\partial$ for any $j\neq i$.
\item if $\epsilon_g(i) = \partial$, then $\epsilon_{g_i^\flat}(i) = \partial$ and $\epsilon_{g_i^\flat}(j) = \epsilon_g(j)$ for any $j\neq i$.
\end{enumerate}
The algebraic semantics of $\mathcal{L}^t := \mathcal{L}(\mathcal{f}^t, \mathcal{G}^t)$ is given by the class of `tense' $\mathcal{L}^t$-algebras, defined as those $\mathcal{L}^t$-algebras $\bba = (L, (\mathcal{F}^t)^{\A}, (\mathcal{G}^t)^{\A})$ such that
\begin{enumerate}
\item for every $f\in \mathcal{F}$ s.t.\ $n_f\geq 1$, all $a_1,\ldots,a_{n_f},b\in L$, and each $1\leq i\leq n_f$,
\begin{itemize}
\item if $\epsilon_f(i) = 1$, then $f(a_1,\ldots,a_i,\ldots a_{n_f})\leq b$ iff $a_i\leq f^\sharp_i(a_1,\ldots,b,\ldots,a_{n_f})$;
\item if $\epsilon_f(i) = \partial$, then $f(a_1,\ldots,a_i,\ldots a_{n_f})\leq b$ iff $a_i\leq^\partial f^\sharp_i(a_1,\ldots,b,\ldots,a_{n_f})$.
\end{itemize}
\item for every $g\in \mathcal{G}$ s.t.\ $n_g\geq 1$, any $a_1,\ldots,a_{n_g},b\in L$, and each $1\leq i\leq n_g$,
\begin{itemize}
\item if $\epsilon_g(i) = 1$, then $b\leq g(a_1,\ldots,a_i,\ldots a_{n_g})$ iff $g^\flat_i(a_1,\ldots,b,\ldots,a_{n_g})\leq a_i$.
\item if $\epsilon_g(i) = \partial$, then $b\leq g(a_1,\ldots,a_i,\ldots a_{n_g})$ iff $g^\flat_i(a_1,\ldots,b,\ldots,a_{n_g})\leq^\partial a_i$.
\end{itemize}
\end{enumerate}

\subsection{Canonical extensions, constructively}
Canonical extensions provide a purely algebraic encoding of Stone-type dualities, and indeed, the existence of the canonical extensions of the best-known varieties of LEs can be proven via preexisting dualities. However, alternative, purely algebraic constructions are available, such as those of \cite{GeHa01,DunnGP05}. These constructions are in fact more general, in that their definition does not rely on principles such as Zorn's lemma. In what follows we will recall them relative to the setting of LEs introduced above.
\begin{definition}\label{Def:Canon:Ext}
Let $\bba$ be a (bounded) sublattice of a complete lattice $\bba'$.
\begin{enumerate}
\item $\bba$ is {\em dense} in $\bba'$ if every element of $\bba'$ can be expressed both as a join of meets and as a meet of joins of elements from $\bba$.
\item $\bba$ is {\em compact} in $\bba'$ if, for all $S, T \subseteq \bba'$, if $\bigvee S\leq \bigwedge T$ then $\bigvee S'\leq \bigwedge T'$ for some finite $S'\subseteq S$ and $T'\subseteq T$.
\item The {\em canonical extension} of a lattice $\bba$ is a complete lattice $\bbas$ containing $\bba$ as a dense and compact sublattice.
\end{enumerate}
\end{definition}
 Let $\kbbas$ and $\obbas$ denote the meet-closure and the join-closure of $\bba$ in $\bbas$ respectively. The elements of $\kbbas$ are referred to as {\em closed} elements, and elements of $\obbas$ as {\em open} elements.
The canonical extension of a bounded lattice $\bba$ exists and is unique up to any isomorphism fixing $\bba$ (cf.\  \cite[Propositions 2.6 and 2.7]{GeHa01}).
In meta-theoretic settings in which Zorn's lemma is available, the canonical extension of a lattice $\bba$ is a {\em perfect} lattice. That is, in addition to being complete, is both completely join-generated by the set $\jira$ of the completely join-irreducible elements of $\bba$, and completely meet-generated by the set $\mira$ of the completely meet-irreducible elements of $\bba$. In our present, constructive setting, canonical extensions might  not be perfect. 

The canonical extension of an LE $\bba$ will be defined as a suitable expansion of the canonical extension of the underlying lattice of $\bba$. Before turning to this definition, recall that taking order-duals interchanges closed and open elements:\label{Comment:Compatibility}
$K({(\bbas)}^\partial) = O(\bbas)$ and $O({(\bbas)}^\partial) =\kbbas$; similarly, $K({(\bba^n)}^\delta) =\kbbas^n$, and $O({(\bba^n)}^\delta) =\obbas^n$. Hence, $K({(\bbas)}^\epsilon) =\prod_i K(\bbas)^{\epsilon(i)}$ and $O({(\bbas)}^\epsilon) =\prod_i O(\bbas)^{\epsilon(i)}$ for every LE $\bba$ and every order-type $\epsilon$ on any $n\in \mathbb{N}$, where
\begin{center}
\begin{tabular}{cc}
$K(\bbas)^{\epsilon(i)}: =\begin{cases}
K(\bbas) & \mbox{if } \epsilon(i) = 1\\
O(\bbas) & \mbox{if } \epsilon(i) = \partial\\
\end{cases}
$ &
$O(\bbas)^{\epsilon(i)}: =\begin{cases}
O(\bbas) & \mbox{if } \epsilon(i) = 1\\
K(\bbas) & \mbox{if } \epsilon(i) = \partial.\\
\end{cases}
$\\
\end{tabular}
\end{center}

As a consequence of these observations, taking the canonical extension of a lattice commutes with taking order duals and products, namely: ${(\bba^\partial)}^\delta = {(\bbas)}^\partial$ and ${(\bba_1\times \bba_2)}^\delta = \bba_1^\delta\times \bba_2^\delta$ (cf.\ \cite[Theorem 2.8]{DunnGP05}). Hence, ${(\bba^\partial)}^\delta$ can be identified with ${(\bbas)}^\partial$, ${(\bba^n)}^\delta$ with ${(\bbas)}^n$, and ${(\bba^\varepsilon)}^\delta$ with ${(\bbas)}^\varepsilon$ for any order-type $\varepsilon$. Thanks to these identifications, in order to extend operations of any arity and which are monotone or antitone in each coordinate from a lattice $\bba$ to its canonical extension, treating the case of {\em monotone} and {\em unary} operations suffices:
\begin{definition}
For every unary, order-preserving operation $f : \bba \to \bba$, the $\sigma$- and $\pi$-{\em extension} of $f$ are defined as follows:
\begin{center}
\begin{tabular}{c}
$f^\sigma (u) =\bigvee \{ \bigwedge \{f(a): k\leq a\in \mathbb{A}\}: u \geq k \in K (\ca)\}$\\
$f^\pi (u) =\bigwedge \{ \bigvee \{f(a): o\geq a\in \mathbb{A}\}: u \leq o \in O (\ca)\}.$\\
\end{tabular}
\end{center}
\end{definition}
It is easy to see that the $\sigma$- and $\pi$-extensions of $\varepsilon$-monotone maps are $\varepsilon$-monotone. Moreover,
the $\sigma$-extension of a map which sends finite (resp.\ finite nonempty) joins or meets in the domain to finite (resp.\ finite nonempty) joins in the
codomain sends {\em arbitrary} (resp.\ {\em arbitrary} nonempty) joins or meets in the domain to {\em arbitrary} (resp.\ {\em arbitrary} nonempty) joins in the codomain. Dually, the $\pi$-extension of a map which sends finite (resp.\ finite nonempty) joins or meets in the domain to finite (resp.\ finite nonempty) meets in the codomain sends {\em arbitrary} (resp.\ {\em arbitrary} nonempty) joins or meets in the domain to {\em arbitrary} (resp.\ {\em arbitrary} nonempty) meets in the codomain (cf.\ \cite[Lemma 4.6]{GeHa01}; notice that the proof given there holds in a constructive meta-theory). Therefore, depending on the properties of the original operation, it is more convenient to use one or the other extension. This justifies the following:
\begin{definition}
The canonical extension of an
$\mathcal{L}$-algebra $\bbA = (L, \mathcal{F}^\bbA, \mathcal{G}^\bbA)$ is the $\mathcal{L}$-algebra $\bbA^\delta: = (L^\delta, \mathcal{F}^{\bbA^\delta}, \mathcal{G}^{\bbA^\delta})$ such that $f^{\bbA^\delta}$ and $g^{\bbA^\delta}$ are defined as the $\sigma$-extension of $f^{\bbA}$ and as the $\pi$-extension of $g^{\bbA}$ respectively, for all $f\in \mathcal{F}$ and $g\in \mathcal{G}$.
\end{definition}
The canonical extension of an LE $\bba$ is a {\em quasi-perfect} LE:
\begin{definition}
\label{def:perfect LE}
An LE $\bbA = (L, \mathcal{F}^\bbA, \mathcal{G}^\bbA)$ is quasi-perfect if $L$ is a complete lattice, and the infinitary versions of the distribution laws defining normality are satisfied for each $f\in \mathcal{F}$ and $g\in \mathcal{G}$.
\end{definition}
\subsection{Adding the fixed point operators}

 In this subsection, we describe two ways of extending  any LE-language by adding fixed point operators. The distinction between the two extensions will become clear when we define their interpretations on LEs. In what follows, $\mathsf{FVAR}$, $\mathsf{PHVAR}$ are disjoint sets of \emph{fixed point variables} and \emph{placeholder variables}, respectively.

For any LE-language $\mathcal{L}$, let $\Tmone$  be the set of terms which extends $\Tm$ by allowing terms $\mu x. t(x)$ and $\nu x.t(x)$ where $t \in \Tmone$, $x \in \mathsf{FVAR}$ and $t(x)$ is positive in $x$. The second extension is denoted $\Tmtwo$ and extends $\Tm$ by allowing construction of the terms $\mu_2 x. t(x)$ and $\nu_2 x.t(x)$ where $t \in \Tmtwo$, $x \in \mathsf{FVAR}$ and $t(x)$ is positive in $x$.

Terms in $\Tm$ are interpreted in LEs as described above. If  $t(x_1,x_2, \ldots,x_n) \in \Tmone$ and $a_1,\ldots,a_{n-1} \in \A$, then
\begin{center}
$\mu x. t(x,a_1,\ldots,a_{n-1}) := \bigwedge \{\,a \in \A \mid t(a,a_1,\ldots,a_{n-1}) \le a \,\}$
\end{center}
if this meet exists, otherwise $\mu x. t(x,a_1,\ldots,a_{n-1})$ is undefined. Similarly,
\begin{center}
$\nu x.t(x,a_1,\ldots,a_{n-1}) := \bigvee \{\, a \in \A \mid a \le t(a,a_1,\ldots,a_{n-1}) \,\}$\end{center}
if this join exists, otherwise $\nu x.t(x,a_1,\ldots,a_{n-1})$ is undefined. For each ordinal $\alpha$ we define $t^\alpha(\bot,a_2,\ldots,a_n)$ as follows:
\begin{center}
\begin{tabular}{l}
$t^0(\bot,a_1,\ldots,a_{n-1}) = \bot, \qquad t^{\alpha+1}(\bot,a_1,\ldots,a_{n-1}) = t\big( t^{\alpha}(\bot,a_1,\ldots,a_{n-1}),a_1,\ldots,a_{n-1}\big)$,\\
$t^{\lambda}(\bot,a_1,\ldots,a_{n-1}) = \bigvee_{\alpha < \lambda} t^{\alpha}(\bot, a_1, \ldots,a_{n-1}) \qquad \text{for limit ordinals } \lambda$;\\
$t_0(\top,a_1,\ldots,a_{n-1}) = \top, \qquad t_{\alpha+1}(\top,a_1,\ldots,a_{n-1}) = t\big( t_{\alpha}(\top,a_1,\ldots,a_{n-1}),a_1,\ldots,a_{n-1}\big)$,\\
$t_{\lambda}(\top,a_1,\ldots,a_{n-1})= \bigwedge_{\alpha < \lambda} t_{\alpha}(\top, a_1, \ldots,a_{n-1}) \qquad \text{for limit ordinals } \lambda$.\\
\end{tabular}
\end{center}

If $t(x_1,\ldots,x_n) \in \Tmtwo$,  then we let
\begin{center}$\LFPtwo x.t(x,a_1,\ldots,a_{n-1}) := \bigvee_{\alpha \geq 0} t^\alpha(\bot,a_1,\ldots,a_{n-1})$\end{center}
\begin{center}$\GFPtwo x.t(x,a_1,\ldots,a_{n-1}):= \bigwedge_{\alpha \geq 0} t_\alpha(\top,a_1,\ldots,a_{n-1})$\end{center}
if this join and  meet exist, otherwise are undefined.

A lattice expansion $\A$ is  \emph{of the first kind} (\emph{of the second kind}) if $t^\A(a_1,\ldots,a_n)$ is defined for all $a_1,\ldots,a_n \in \A$ and all $t \in \Tmone$ ($t \in \Tmtwo$). Henceforth we will refer to these algebras as \emph{mu-algebras of the first kind} (\emph{of the second kind}).
When restricted to the Boolean case, our mu-algebras of the first kind are essentially the modal mu-algebras defined in~\cite[Definition 2.2]{BH12} and~\cite[Definition 5.1]{AmKwMe95}. 
Every  mu-algebra of the second kind is a mu-algebra of the first kind.
(cf.\  \cite[Proposition 2.4]{AmKwMe95} and \cite[Lemma 2.2]{CoCr14}. These proofs straightforwardly extend to the setting of general LEs).
 Hence, the interpretation of the two types of fixed point binders  on mu-algebras of the second kind will agree. That is, $\mu X.\phi(X) =\mu_2 X.\phi(X)$ and $\nu X.\psi(X) = \nu_2 X.\psi(X)$ in mu-algebras of the second kind.


The final sets of terms, $\Tm_*$ (resp.\ $\Tm^t_*$), are obtained as extensions of $\mathcal{L}$ (resp.\ $\mathcal{L}^t$) by allowing $\mu^* x. s(x)$ and $\nu^* x.s(x)$ whenever $s \in \Tm_*$ (resp.\ $s \in \Tm^t_*$) and is positive in $x$. Terms in $\Tm_*$ and $\Tm^t_*$ are only interpreted in the constructive canonical extensions $\A^\delta$ of lattice expansions $\A$. If $s(x_1,x_2, \ldots,x_n) \in \Tm_* \cup \Tm^t_*$ and $a_1,\ldots,a_{n-1} \in \A^\delta$, then $\mu^* x_1. s(x_1,a_1,\ldots,a_{n-1}) := \bigwedge \{\,a \in \A \mid s(a,a_1,\ldots,a_{n-1}) \le a \,\}$ and $\nu^* x_1.s(x_1,a_2,\ldots,a_n) := \bigvee \{\, a \in \A \mid a \le s(a,a_1,\ldots,a_{n-1}) \,\}$.
As the canonical extension $\A^\delta$ is a complete lattice, the interpretation of $\mu^* x. t(x)$ or $\nu^* x.t(x)$ is always defined.
 For any term $\phi \in \Tm^t_1$ we let $\phi^*$ denote the $\Tm^t_*$ term obtained from $\phi$ by replacing all occurrences of $\mu$ and $\nu$ with $\mu^*$ and $\nu^*$, respectively.
The main feature of the $\mu^*$ and $\nu^*$ binders  is  that their interpretation does not change  from $\A$ to $\A^\delta$.

\subsection{The language of constructive ALBA for LEs}\label{Subsec:Expanded:Land}

The expanded language manipulated by the algorithm $\mu^*$-ALBA (cf.\ section \ref{sec:constr:mu}) includes the $\mathcal{L}^t$-connectives, as well as a denumerably infinite set of sorted variables $\mathsf{NOM}$ called {\em nominals}, a denumerably infinite set of sorted variables $\mathsf{CO\text{-}NOM}$, called {\em co-nominals}, and fixed point binders (depending on the language). We let $\nomi, \nomj$  denote nominals,  and $\cnomm, \cnomn$ denote co-nominals. While in the non-constructive setting nominals and co-nominals range over the completely join-irreducible and the completely meet-irreducible elements of perfect LEs, respectively, in the present, constructive setting, nominals and co-nominals will be interpreted as elements of $\kbbas$ and $\obbas$, respectively.

Formulas in the extended language $\mathcal{L}_1^{+}$ are defined by the following recursion:
\begin{center}
$\phi ::= \bot \mid \top \mid p \mid X \mid \mathbf{j} \mid \mathbf{m} \mid \phi \wedge \psi \mid \phi \vee \psi \mid f(\overline{\phi}) \mid g(\overline{\phi}) \mid \mu X.\phi(X)\mid \nu X.\phi(X)$
\end{center}
where $p \in \mathsf{PROP}$, $X \in \mathsf{FVAR}$, $\nomj \in \mathsf{NOM}$, $\cnomm \in \CNOM$, $f\in\mathcal{F}^{t}$, $g\in\mathcal{G}^{t}$, and  $\phi$ is positive in $X$ in $\mu X.\phi(X)$ and $\nu X.\phi(X)$. Formulas in $\mathcal{L}_2^{+}$ (resp.\ $\mathcal{L}_*^{+}$) are defined by replacing the fixed point operators $\mu,\nu$ with $\mu_2,\nu_2$ (resp.\ $\mu^*, \nu^*$).

Placeholder variables from $\mathsf{PHVAR}$, denoted $x, y, z$, will be used as generic variables which can take on the roles of propositional and fixed point variables. They will also be used to enhance the clarity of the exposition when dealing with substitution instances of formulas. Let $\tau$ be an order-type over $n$.
An $\mathcal{L}_1^{+}$ (resp.\ $\mathcal{L}_2^{+}$, $\mathcal{L}_*^{+}$) formula is \emph{pure} if it contains no ordinary (propositional) variables but only, possibly, nominals and co-nominals, and is an $\mathcal{L}_1^{+}$ (resp.\ $\mathcal{L}_2^{+}$, $\mathcal{L}_*^{+}$)-\emph{sentence} if it contains no free fixed point variables.

A \emph{quasi-inequality} of $\mathcal{L}_1^{+}$ (resp.\ $\mathcal{L}_2^{+}$, $\mathcal{L}_*^{+}$) is an expression of the form $\phi_1 \leq \psi_1 \amp \cdots \amp \phi_n \leq \psi_n \Rightarrow \phi \leq \psi$ where the $\phi_i$, $\psi_i$, $\phi$ and $\psi$ are formulas of $\mathcal{L}_1^{+}$ (resp.\ $\mathcal{L}_2^{+}$, $\mathcal{L}_*^{+}$).

\paragraph{Quasi-inequalities, assignments, validity.} Constructive canonical extensions of ${\mathcal{L}}$-algebras can be naturally endowed with the structure of ${\mathcal{L}^+}$-algebras (cf.\ Section 2.6 \cite{CoPaZh16b}). Building on this fact, we can use constructive canonical extensions of ${\mathcal{L}}$-algebras as a semantic environment for the language ${\mathcal{L}^+}$ as follows.
 For any $\mathcal{L}$-algebra $\bba$, an \emph{assignment} on $\bba$ sends propositional variables to elements of $\bba$ and is extended to formulas of $\mathcal{L}_1$, $\mathcal{L}_2$ and $\mathcal{L}_*$ in the usual way, where these are defined. An  \emph{assignment} on $\bbas$ is a map $v: \mathsf{PROP} \cup \mathsf{NOM} \cup \mathsf{CO\mbox{-}NOM} \rightarrow \bbas$ sending propositional variables to elements of $\bbas$, nominals to $\kbbas$ and co-nominals to $\obbas$ and extends to all formulas of $\mathcal{L}^+_*$, $\mathcal{L}^+_1$ and $\mathcal{L}^+_2$. An \emph{admissible assignment} on $\bbas$ \label{admissible:assignment} is an assignment which takes all propositional variables to elements of $\bba$. An $\mathcal{L}^+$-inequality $\alpha \leq \beta$ is \emph{admissibly valid} on $\bba$, denoted $\bbas \models_\bba \alpha \leq \beta$, if it holds under all admissible assignments.
 A quasi-inequality $\phi_1 \leq \psi_1 \amp \cdots \amp \phi_1 \leq \psi_1 \Rightarrow \phi \leq \psi$ is satisfied under an assignment $V$ in an algebra $\A$ of the appropriate sort, written $\A, V \models \phi_1 \leq \psi_1 \amp \cdots \amp \phi_n \leq \psi_n \Rightarrow \phi \leq \psi$ if $\A, V \not \models \phi_i \leq \psi_i$ for some $1 \leq i \leq n$ or $\A, v \models \phi \leq \psi$. A quasi-inequality is (admissibly) valid in an algebra if it is satisfied by every (admissible) assignment.

\paragraph{Signed generation trees.}
For any formula/term $\phi$ in $\mathcal{L}_1^+$ and $\mathcal{L}_*^+$, we assign two \emph{signed generation trees}
$+\phi$ and $-\phi$. The generation tree is
constructed as usual, beginning at the root with the main connective and then
branching out into $n$-nodes at each $n$-ary connective. Each leaf is either a propositional
variable, a fixed point variable, or a constant. Each node is signed as follows:
\begin{itemize}
\item the root node of $+\phi$ is signed $+$ and the root node of $-\phi$ is signed $-$;
\item if a node is $\vee, \wedge$ assign the same sign to its successor nodes;
\item if a node is $h\in \mathcal{F}^t\cup\mathcal{G}^t$, assign the same (resp.\ the opposite) sign to every node corresponding to a coordinate $i$ such that $\epsilon_h(i)= 1$ (resp.\ $\epsilon_h(i)= \partial$);
\item if a node is $\mu x.\varphi(x)$, $\mu^* x.\varphi(x)$, $\nu x.\varphi(x)$ or $\nu^* x.\varphi(x)$ (with every free occurrence of $x$ in the positive generation tree of $\varphi$ labelled positively)
then assign the same sign to the successor node.
\end{itemize}
A node in a signed generation tree is
\emph{positive} if it is signed ``$+$'' and \emph{negative} if it is signed ``$-$''. 

\section{Two kinds of canonicity, constructively}\label{Sec:TwoKindsCanon}

In this section we give a brief conceptual and methodological overview of the main results of this paper.

The arguments for canonicity have the ``U-shaped'' format typical of the unified correspondence paradigm (see \cite{CoGhPa14}) as generically illustrated in Figure~\ref{fig:generalU}.
\begin{figure}
\begin{center}
\begin{tikzpicture}
\path (0,0) node(a) {$\A \models \alpha \le \beta$}
			(0,-0.5) node(b) [rotate=90] {$\Leftrightarrow$}
			(0,-1) node(c) {$\A^{\delta} \models_\A \alpha \le \beta$}
			(0,-1.5) node(i) [rotate=90] {$\Leftrightarrow$}
			(0,-2) node(d) {$\A^{\delta} \models_\A \mathrm{ALBA}(\alpha \le \beta)$}
			(5,-2) node(f) {$\A^\delta \models \mathrm{ALBA}(\alpha \le \beta)$}
			(5,-1) node(g) [rotate=90] {$\Longleftrightarrow$}
			(2.5,-2) node(j) {$\Longleftrightarrow$}
			(5,0) node(k) {$\A^\delta \models \alpha \le \beta$};
\end{tikzpicture}
\caption{The U-shaped argument for canonicity of inequalities interpreted on a LE $\A$.}
\label{fig:generalU}
\end{center}
\end{figure}
Going down the left-hand arm of the diagram, the first bi-implication is justified by the fact that validity in $\A$ is the same as admissible validity in $\A^\delta$, provided we make some stipulations regarding  the interpretation of fixed point binders. The extended language $\mathcal{L}^+_*$ can be interpreted in $\A^\delta$ with nominals and co-nominals running over closed and open elements, respectively.

The inequality is now equivalently transformed into a set of pure (quasi-) in\-e\-qua\-li\-ties, denoted $\mathrm{ALBA}(\alpha \le \beta)$ in Figure~\ref{fig:generalU}. This is done by means of a calculus of rewrite rules encapsulated in the constructive $\mu^*$-ALBA, which we present in Section~\ref{sec:constr:mu}. The fact that admissible and ordinary validity coincide for pure inequalities allows us to traverse the the bi-implication forming the base of the ``U''.

We proceed up the right-hand arm of the ``U'' by reversing the rewrite rules applied when going down the left-hand side. The equivalences are justified by the fact that these rules preserve validity on quasi-perfect algebras (cf.\ Definition \ref{def:perfect LE}).

If $\phi$ is a formula without fixed point binders, then the term function $\phi^{\A^{\delta}}$ \emph{extends} the term function $\phi^{\A}$, i.e., they agree on arguments from $\bba$. This  is usually of crucial importance in proving that an equation is canonical. As soon as we add fixed point binders, this extension property fails. Indeed, $(\phi(X))^{\A^{\delta}}$ can have more pre-fixed points in $\A^{\delta}$ than $(\phi(X))^{\A}$ has in $\A$, and so $(\mu X. \phi(X))^{\A^{\delta}}$ would generally be smaller than $(\mu X. \phi(X))^{\A}$. This phenomenon creates significant obstacles for standard canonicity arguments, and is not an additional difficulty posed by the constructive environment. 

One way around these difficulties, adopted by Bezhanishvili and Hodkinson~\cite{BH12},  is to require that only pre-fixed points from the smaller algebra $\A$ are used in calculating $(\mu X. \phi(X))^{\A^{\delta}}$. This is tantamount to interpreting $\mu X. \phi(X)$ in $\A^{\delta}$ as $(\mu^* X. \phi(X))^{\A^{\delta}}$.
Thus, in \cite{CoCr14},  two different notions of canonicity for the mu-calculus were considered, and their two ensuing canonicity results were shown for certain classes of mu-inequalities. The counterparts of these results hold in a constructive general lattice environment. Specifically,
following \cite{CoCr14}, we call an inequality $\phi \le \psi$ \emph{tame canonical} when $\A \models \phi \le \psi$ if and only if $\A^\delta \models \varphi^* \le \psi^*$ for all $\mu$-algebras $\A$ of the first kind.  We generalize the \emph{tame inductive} mu-inequalities of \cite{CoCr14} to the LE-setting, and prove that they are tame canonical in a constructive meta-theory.

Of course, the usual notion of canonicity may also be applied to formulas with fixed-point binders, i.e., that $\A \models \phi \le \psi$ iff $\A^\delta \models \phi \le \psi$, where fixed points are interpreted in the standard way, e.g.\ least fixed points in $\A^\delta$ are calculated as the meet of all pre-fixed points \emph{in $\A^\delta$}. A canonicity result of this kind can be proved by generalizing the class of \emph{restricted inductive mu-inequalities} of \cite{CoCr14} to the LE-setting, and showing that they are preserved under constructive canonical extensions of mu-algebras of the second kind.

Whenever a \emph{tame} run of $\mu^*$-ALBA succeeds on a mu-inequality $\phi \le \psi$, we have that $\phi \le \psi$ is \emph{tame canonical}.
Moreover, whenever a \emph{proper} run succeeds on a mu-inequality $\alpha \le \beta$, then $\alpha \le \beta$ will be \emph{canonical}.
Finally, for every tame inductive mu-inequality (respectively, a restricted inductive mu-inequality), there exists a tame (respectively, proper) run of $\mu^*$-ALBA which succeeds on that inequality.

\section{Syntactic classes}
\label{sec:synclass}

In this section, we introduce some syntactically defined classes of mu-inequalities, the most general of which is the counterpart, in the language of normal LEs plus fixed points, of the recursive mu-inequalities introduced in \cite{CoFoPaSo15}. In a constructive setting, the members of this class will all have correspondents in a first-order language plus fixed points. The remaining two classes of mu-inequalities defined in this section are subclasses of the class of recursive mu-inequalities, and are those for which the two canonicity results hold.


For any $\Tmone$-sentence $\varphi(p_1,\ldots p_n)$, any order-type $\epsilon$ over $n$, and any $1\leq i\leq n$, an \emph{$\epsilon$-critical node} in a signed generation tree of $\varphi$ is a (leaf) node $+p_i$ with $\epsilon_i = 1$, or $-p_i$ with $\epsilon_i = \partial$. An $\epsilon$-{\em critical branch} in the tree is a branch terminating in an $\epsilon$-critical node. The intuition, which will be built upon later, is that variable occurrences corresponding to $\epsilon$-critical nodes are \emph{to be solved for, according to $\epsilon$}.

In the signed generation tree of a $\Tmone$-sentence $\varphi(p_1,\ldots p_n)$ a \emph{live branch} is a branch ending in a (signed) propositional variable. In particular, all critical branches are live. A branch is not live iff it ends in a propositional constant ($\top$ or $\bot$) or in a fixed point variable.

For every $\Tmone$-sentence $\varphi(p_1,\ldots p_n)$, and every order-type $\epsilon$, we say that $+\varphi$ (resp. $-\varphi$) {\em agrees with} $\epsilon$, and write $\epsilon(+\varphi)$ (resp.\ $\epsilon(-\varphi)$), if every leaf node in the signed generation tree of $+\varphi$ (resp.\ $-\varphi$) which is labelled with a propositional variable is $\epsilon$-critical.
We will also make use of the {\em sub-tree relation} $\gamma\prec \varphi$, which extends to signed generation trees, and we will write $\epsilon(\gamma)\prec \ast \varphi$ to indicate that $\gamma$, regarded as a sub- (signed generation) tree of $\ast \varphi$, agrees with $\epsilon$.

\begin{table}
\begin{center}
\begin{tabular}{| c | c | c | c |}
\hline
Outer Skeleton ($P_3$) &Inner Skeleton ($P_2$) &PIA ($P_1$)\\
\hline
$\Delta$-adjoints &Binders &Binders\\
\begin{tabular}{ c c c c  }
&$+$ &$\vee$ &  \\
&$-$ &$\wedge$ &\\
\hline
\end{tabular}
&\begin{tabular}{c c c c ccccc}
& &&$+$ &$\mu$ &&\\
& &&$-$ &$\nu$ &&\\
\hline
\end{tabular}
&\begin{tabular}{c c c c c c c c }
& & &$+$ &$\nu$ & & \\
& & &$-$ &$\mu$ & & \\
\hline
\end{tabular}
\\
SLR &SLA &SRA\\
\begin{tabular}{c c c c}
& $+$ &$f$ &\\
&$-$ &$g$ &\\
\end{tabular}
&\begin{tabular}{c c c c c c c c}
&$+$ &$\vee$ &$f$ &$(n_f = 1)$ &\\
&$-$ &$\wedge$ &$g$ &$(n_g = 1)$ &\\
\hline
\end{tabular}
&\begin{tabular}{c c c c c c c c}
&$+$ &$\wedge$ &$g$ &$(n_g = 1)$ &\\
&$-$ &$\vee$ &$f$ &$(n_f = 1)$ &\\
\hline
\end{tabular}
\\
&SLR &SRR\\
&\begin{tabular}{c c c c c c c c}
$+$ &$\phantom{\wedge}$ &$f$ &$(n_f\geq 2)$ &\\
$-$ & &$g$ &$(n_g\geq 2)$ &\\
\end{tabular}
&\begin{tabular}{c c c c c cc}
$+$ &$\phantom{\wedge}$&$g$ &$(n_g\geq 2)$ &\\
$-$ &&$f$ &$(n_f\geq 2)$ &\\
\end{tabular}
\\
\hline
\end{tabular}
\end{center}
\caption{Skeleton and PIA nodes.}\label{Skeleton:PIA:Node:Table}
\end{table}


\begin{definition}
Nodes in signed generation trees will be called \emph{skeleton nodes} and \emph{PIA nodes} and further classified as $\Delta$-adjoint, SLR, Binders, SLA, SRA or SRR, according to the specification given in Table \ref{Skeleton:PIA:Node:Table}.\footnote{The abbreviations SLR, SLA, SRA and SRR stand for syntactically left residual, left adjoint, right adjoint and right residual, respectively. Nodes are thus classified according to the order-theoretic properties of their interpretations, see \cite{CoGhPa14,ConPal13} for further discussion on methodology and nomenclature.}
\label{def:good branch}
Let $\phi(p_1,\ldots,p_n)$ be a formula in the propositional variables $p_1, \ldots, p_n$,
and let
$\epsilon$ be an order-type on $\{1, \ldots, n\}$.

A branch in a signed generation tree $\ast \varphi$, for $\ast \in \{+, - \}$, ending in a propositional variable is an $\epsilon$-\emph{good branch} if, apart from the leaf, it is the concatenation of three paths $P_1$, $P_2$, and $P_3$, each of which may possibly be of length $0$, such that $P_1$ is a path from the leaf consisting only of PIA-nodes, $P_2$ consists only of inner skeleton-nodes, and $P_3$ consists only of outer skeleton-nodes and, moreover, it satisfies conditions (GB1), (GB2) and (GB3), below.

\begin{description}

\item[{\bf (GB1)}] The formula corresponding to the uppermost node on $P_1$ is a sentence.
\item[{\bf (GB2)}] For every SRR-node in $P_1$ of the form $h(\overline{\gamma}, \beta)$, where $\beta$ is the coordinate where the branch lies, every $\gamma$ in $\overline{\gamma}$ is a mu-sentence and $\epsilon^\partial (\gamma)\prec \ast\varphi$ (i.e., each $\gamma$ contains no variable occurrences to be solved for --- see above).
\item[{\bf (GB3)}] For every SLR-node in $P_2$ of the form $h(\overline{\gamma}, \beta)$, where $\beta$ is the coordinate where the branch lies, every $\gamma$ in $\overline{\gamma}$ is a mu-sentence and $\epsilon^\partial (\gamma)\prec \ast\varphi$ (see above for this notation).\\
\end{description}
\end{definition}
%
%
\noindent Our main interest is in
$\epsilon$-good branches satisfying some of the
additional properties in the following definition.

\begin{definition}\label{def:nbpia-omconf}
Let $\phi(p_1,\ldots,p_n)$ be a formula, 
$\epsilon$ be an order-type on $\{1, \ldots, n\}$, and $<_{\Omega}$ a strict partial order on the variables $p_1,\ldots p_n$. An $\epsilon$-good branch may satisfy one or more of the following properties:

\begin{description}
\item[{\bf (NB-PIA)}] $P_1$ contains no fixed point binders.


\item[{\bf (NL)}] For every SLR-node in $P_2$ of the form $h(\overline{\gamma}, \beta)$, where $\beta$ is the coordinate where the branch lies, the signed generation tree of each $\gamma$ contains no live branches.
\item[{\bf ($\Omega$-CONF)}] For every SRR-node in $P_1$ of the form $h(\overline{\gamma}, \beta)$, where $\beta$ is the coordinate where the branch lies, $p_j <_{\Omega} p_i$ for every $p_j$ occurring in $\gamma$, where $p_i$ is the propositional variable labelling the leaf of the branch.
\end{description}
\end{definition}

\begin{definition}\label{def:rstrctd:tm:indctv}
 For any order-type $\epsilon$ and  strict partial order $<_{\Omega}$ on the variables $p_1,\ldots p_n$, the signed generation tree $\ast \phi$, $\ast \in \{-, + \}$, of a term $\phi(p_1,\ldots p_n)$ is called
\begin{enumerate}
\item \emph{$\epsilon$-recursive} if every $\epsilon$-critical branch is $\epsilon$-good.
\item \emph{$(\Omega, \epsilon)$-inductive} it is $\epsilon$-recursive and every $\epsilon$-critical branch satisfies ($\Omega$-CONF).
\item \emph{restricted $(\Omega,\epsilon)$-inductive} if it is $(\Omega,\epsilon)$-inductive and
 \begin{enumerate}
 \item every $\epsilon$-critical branch satisfies (NB-PIA) and (NL),
 \item every occurrence of a binder is on an $\epsilon$-critical branch.
 \end{enumerate}
\item \emph{tame $(\Omega,\epsilon)$-inductive} if it is $(\Omega,\epsilon)$-inductive and
 \begin{enumerate}
 \item $\Omega = \varnothing$,
 \item no binder occurs on any $\epsilon$-critical branch,
 \item the only nodes involving binders which are allowed to occur are $+\nu$ and $-\mu$.
 \end{enumerate}
\end{enumerate}
An inequality $\phi \leq \psi$ is $\epsilon$-recursive (resp., $(\Omega, \epsilon)$-inductive, restricted $(\Omega,\epsilon)$-inductive, tame $(\Omega,\epsilon)$-inductive) if $+\phi$ and $-\psi$ are both $\epsilon$-recursive (resp., $(\Omega, \epsilon)$-inductive, restricted $(\Omega,\epsilon)$-inductive, tame $(\Omega,\epsilon)$-inductive).

An inequality $\phi \leq \psi$ is recursive (resp., inductive, restricted inductive, tame inductive) if $\phi \leq \psi$ is $\epsilon$-recursive (resp., $(\Omega, \epsilon)$-inductive, restricted $(\Omega,\epsilon)$-inductive, tame $(\Omega,\epsilon)$-inductive) for some strict partial order $\Omega$ and order-type $\epsilon$.

The corresponding classes of inequalities will be referred to as the \emph{recursive} (resp., \emph{inductive, restricted inductive, tame inductive})
\emph{mu-inequalities}, or the \emph{recursive} (resp., \emph{inductive, restricted inductive, tame inductive}) \emph{mu-formulas}, if the inequality signs have been replaced with implications.
\end{definition}

\section{Constructive $\mu^{*}$-ALBA}\label{sec:constr:mu}
In this section, we introduce the constructive and general lattice version of the algorithm $\mu^{*}$-ALBA introduced in \cite{CoCr14}, which in its turn is a restricted  version of the algorithm $\mu$-ALBA,  introduced in~\cite{CoFoPaSo15} to calculate first-order correspondents of recursive inequalities  from bi-intuitionistic modal mu-calculus. Constructive $\mu^{*}$-ALBA is the fundamental tool to prove the canonicity results via the argument discussed in Section \ref{Sec:TwoKindsCanon}.  As usual, the goal of constructive $\mu^{*}$-ALBA is to \emph{eliminate propositional variables} from inequalities, while maintaining admissible validity. The purpose of this is to make the transition from \emph{admissible validity} to \emph{validity} in the argument for canonicity. Below, we outline its general strategy.

Constructive $\mu^{*}$-ALBA, from now on abbreviated as ALBA, takes an $\mathcal{L}_1$-inequality $\phi \leq \psi$ as input, and proceeds in three stages. The first stage preprocesses by eliminating all uniformly occurring propositional variables, applying distribution rules for $f\in \mathcal{F}$ and $g\in \mathcal{G}$ and splitting rules exhaustively, and  converting all occurrences of $\mu X.\phi(X)$ to $\mu^* X.\phi(X)$ and all occurrences of $\nu X. \psi(X)$ to $\nu^* X.\psi(X)$. We emphasize that this step is required in both tame \emph{and} proper runs of ALBA (see Appendix \ref{sec: appendix ALBA}). The preprocessing produces a finite set of inequalities, $\phi'_i \leq \psi'_i$, $1 \leq i \leq n$.

Now ALBA forms the \emph{initial quasi-inequalities} $\bigamp S_i \Rightarrow \sf{Ineq}_i$, compactly represented as tuples $(S_i, \sf{Ineq}_i)$ referred to as \emph{systems}, with each $S_i$ initialized to the empty set and $\sf{Ineq}_i$ initialized to $\phi'_i \leq \psi'_i$.

The second stage (called the reduction stage) transforms $S_i$ and $\mathsf{Ineq}_i$ through the application of transformation rules (see Appendix \ref{sec: appendix ALBA}). The aim is to eliminate all propositional variables from $S_i$ and $\mathsf{Ineq}_i$ in favour of nominals and co-nominals. A system for which this has been done will be called \emph{pure} or \emph{purified}. The actual eliminations are effected through the Ackermann rules, while the other rules are used to bring $S_i$ and $\mathsf{Ineq}_i$ into the appropriate shape which make these applications possible. Once all propositional variables have been eliminated, this phase terminates and returns the pure quasi-inequalities $\bigamp S_i \Rightarrow \mathsf{Ineq}_i$.

The third stage either reports failure if some system could not be purified, or else returns the conjunction of the pure quasi-inequalities $\bigamp S_i \Rightarrow \mathsf{Ineq}_i$, which we denote by $\mathsf{ALBA}(\phi \leq \psi)$.

A more complete outline of each of the three stages will be given in Appendix \ref{sec: appendix ALBA}.


\section{Main results}\label{Sec:Main:Results}

In this section we state the main results of the paper, namely that all restricted inductive inequalities are constructively canonical and that all tame inductive inequalities are constrictively tame canonical. The proof strategy is the same in both cases: one first proves, by means of a `U-shaped' argument as discussed in Section \ref{Sec:TwoKindsCanon}, that successful runs of constructive ALBA satisfying certain conditions guarantee these types of canonicity. Then it is shown that all members of the two classes of inequalities are successfully reducible by means of runs respectively satisfying these properties. The main canonicity results then follow as corollaries. Because of space limitations we do not include the proofs.

A proof combining insights from the proofs of \cite[Theorems 9.9 and 9.4]{CoCr14} and \cite[Theorem 7.1]{CoPaZh16b} suffices to establish the following two propositions:

\begin{prop}\label{thm:truecanonicity}
Let $\A$ be a mu-algebra of the second kind and let $\phi \le \psi$ be an $\Tm_1$-inequality
on which a proper and pivotal run of $\mu^*$-ALBA succeeds.
If $\A \models \phi \le \psi$ then $\A^\delta \models \phi \le \psi$.
\end{prop}

\begin{prop}\label{thm:tamecanonicity}
All $\mathcal{L}_1$-inequalities on which a tame and pivotal run of constructive $\mu^*$-ALBA succeeds are constructively tame canonical.
\end{prop}

Again, generalizing to the non-distributive environment and amalgamating the strategies from  \cite[Section 10]{CoCr14} and \cite[Section 6]{CoPaZh16b} yields the next proposition:

\begin{prop}\label{prop:ALBA:suc:rest}
Constructive $\mu^*$-ALBA succeeds on all restricted inductive $\mathcal{L}$-inequalities by means of proper and pivotal runs.
\end{prop}

\begin{prop}\label{prop:ALBA:suc:tame}
Constructive $\mu^*$-ALBA succeeds on all tame inductive $\mathcal{L}$-inequalities by means of tame and pivotal runs.
\end{prop}

The canonicity of restricted inductive $\mathcal{L}$-inequalities now follows as a corollary of Propositions \ref{thm:truecanonicity} and \ref{prop:ALBA:suc:rest}.

\begin{theorem}
All restricted inductive $\mathcal{L}$-inequalities are constructively canonical over mu-algebras of the second kind.
\end{theorem}

Similarly, the tame canonicity of all tame inductive inequalities follows from Propositions \ref{thm:tamecanonicity} and  \ref{prop:ALBA:suc:tame}.

\begin{theorem}
All tame inductive $\mathcal{L}$-inequalities are constructively canonical over mu-algebras of the first kind.
\end{theorem}

\bibliographystyle{aiml16}
\bibliography{aiml16}

\Appendix

\section{Stages and rules of $\mu^{*}$-ALBA}
\label{sec: appendix ALBA}
\subsection{Stage 1: Preprocessing and initialization} ALBA receives an $\mathcal{L}_1$-inequality $\phi \leq \psi$ as input. It applies the following {\bf rules for elimination of monotone variables} to $\phi \leq \psi$ exhaustively, in order to eliminate any propositional variables which occur uniformly:

\begin{prooftree}
\AxiomC{$\alpha(p) \leq \beta(p)$}\UnaryInfC{$\alpha(\bot) \leq \beta(\bot)$}
\AxiomC{$\gamma(p) \leq \delta(p)$}\UnaryInfC{$\gamma(\top) \leq \delta(\top)$}
\noLine\BinaryInfC{}
\end{prooftree}
for $\alpha(p) \leq \beta(p)$ positive and $\gamma(p) \leq \delta(p)$ negative in $p$, respectively.\footnote{\label{footnote:uniformterms}
A term $\phi$ is \emph{positive} (\emph{negative}) in a variable
$p$ if in the generation tree $+\phi$ all $p$-nodes are signed $+$ ($-$). An inequality $\phi \le \psi$ is \emph{positive} (\emph{negative})
in  $p$ if $\phi$ is negative (positive) in $p$ and $\psi$ is positive (negative) in $p$.
}

Next, ALBA exhaustively distributes $f\in \mathcal{F}$ over $+\vee$, and $g\in \mathcal{F}$ over $-\wedge$, so as to bring occurrences of $+\vee$ and $-\wedge$ to the surface wherever this is possible, and then eliminate them via exhaustive applications of {\em splitting} rules.
\paragraph{Splitting-rules.}

\begin{prooftree} \AxiomC{$\alpha \leq \beta \wedge \gamma
$}\UnaryInfC{$\alpha \leq \beta \quad \alpha \leq \gamma$}
\AxiomC{$\alpha \vee \beta \leq \gamma$}\UnaryInfC{$\alpha \leq \gamma \quad \beta \leq \gamma$}
\noLine\BinaryInfC{}
\end{prooftree}

This gives rise to a set of inequalities $\{\phi_i' \leq \psi_i'\mid 1\leq i\leq n\}$. For each of them, ALBA converts all occurrences of $\mu X.\phi(X)$ to $\mu^* X.\phi(X)$ and all occurrences of $\nu X. \psi(X)$ to $\nu^* X.\psi(X)$, and
forms the \emph{initial quasi-inequality} $\bigamp S_i \Rightarrow \sf{Ineq}_i$, compactly represented as a tuple $(S_i, \sf{Ineq}_i)$ referred as \emph{initial system}, with each $S_i$ initialized to the empty set and $\sf{Ineq}_i$ initialized to $\phi'_i \leq \psi'_i$. Each initial system is passed separately to stage 2, described below, where we will suppress indices $i$.

\subsection{Stage 2: Reduction and elimination}

The aim of this stage is to eliminate all occurring propositional variables from a given system $(S, \mathsf{Ineq})$. This is done by means of the following \emph{approximation rules}, \emph{residuation rules}, \emph{splitting rules}, and \emph{Ackermann rules}, collectively called \emph{reduction rules}. The terms and inequalities in this subsection are from $\mathcal{L}_*^{+}$.

\paragraph{Approximation rules.} There are four approximation rules. Each of these rules functions by simplifying $\mathsf{Ineq}$ and adding an inequality to $S$. We write $\alpha(!x)$ to indicate that the placeholder variable $x$ has a unique occurrence in formula $\alpha$.

\begin{description}
\item[Left-positive approximation rule.] $\phantom{a}$
\begin{center}
\AxiomC{$(S, \;\; \phi'(\gamma / !x)\leq \psi)$}
\RightLabel{$(L^+A)$}
\UnaryInfC{$(S\! \cup\! \{ \nomj \leq \gamma\},\;\; \phi'(\nomj / !x)\leq \psi)$}
\DisplayProof
\end{center}
with $+x \prec +\phi'(!x)$, the branch of $+\phi'(!x)$ starting at $+x$ subject to the restrictions detailed below, 
$\gamma$ belonging to the smaller language $\mathcal{L}_*$
%
 %
 %
 and $\nomj$ being the first nominal variable not occurring in $S$ or $\phi'(\gamma / !x)\leq \psi$.
\item[Left-negative approximation rule.]$\phantom{a}$
\begin{center}
\AxiomC{$(S, \;\; \phi'(\gamma / !x)\leq \psi)$}
\RightLabel{$(L^-A)$}
\UnaryInfC{$(S\! \cup\! \{ \gamma\leq\cnomm\},\;\; \phi'(\cnomm / !x)\leq \psi)$}
\DisplayProof
\end{center}
%
%
with $-x \prec +\phi'(!x)$, the branch of $+\phi'(!x)$ starting at $-x$ subject to the restrictions detailed below, 
$\gamma$ belonging to the smaller language $\mathcal{L}_*$ and
 $\cnomm$ being the first co-nominal not occurring in $S$ or $\phi'(\gamma / !x)\leq \psi$.
\item[Right-positive approximation rule.]$\phantom{a}$
\begin{center}
\AxiomC{$(S, \;\; \phi\leq \psi'(\gamma / !x))$}
\RightLabel{$(R^+A)$}
\UnaryInfC{$(S\! \cup\! \{ \nomj \leq \gamma\},\;\; \phi\leq \psi'(\nomj / !x))$}
\DisplayProof
\end{center}
%
%
with $+x \prec -\psi'(!x)$, the branch of $-\psi'(!x)$ starting at $+x$ subject to the restrictions detailed below, 
$\gamma$ belonging to the smaller language $\mathcal{L}_*$ and
 $\nomj$ being the first nominal not occurring in $S$ or $\phi\leq \psi'(\gamma / !x)$.
\item[Right-negative approximation rule.] $\phantom{a}$
\begin{center}
\AxiomC{$(S, \;\; \phi\leq \psi'(\gamma / !x))$}
\RightLabel{$(R^-A)$}
\UnaryInfC{$(S\! \cup\! \{ \gamma\leq \cnomm\},\;\; \phi\leq \psi'(\cnomm / !x))$}
\DisplayProof
\end{center}
with $-x \prec -\psi'(!x)$, the branch of $-\psi'(!x)$ starting at $-x$ subject to the restrictions detailed below, 
$\gamma$ belonging to the smaller language $\mathcal{L}_*$ and $\cnomm$ being the first co-nominal not occurring in $S$ or $\phi\leq \psi'(\gamma / !x))$.
\end{description}

The restrictions on $\phi'$ and $\psi'$ in the approximation rules above are formulated in terms of the following:

\begin{definition}
 For any mu-algebra $\mathbb{C}$ and order-type $\tau$, a join $\bigvee S$ in $(\mathbb{C}^{\delta})^\tau$ is called \emph{$\mathbb{C}^\tau$-targeted} if $\bigvee S \in \mathbb{C}^\tau$. A map $f: (\mathbb{C}^{\delta})^\tau \rightarrow \mathbb{C}^{\delta}$ \emph{preserves $\mathbb{C}^\tau$-targeted joins} if $f(\bigvee S) = \bigvee_{s \in S}f(s)$ for every $S \subseteq (\mathbb{C}^{\delta})^\tau$ such that $\bigvee S$ is $\mathbb{C}^\tau$-targeted.  \emph{Targeted meets} and their preservation are defined order-dually.
\end{definition}
Let us now list the requirements on $\phi'$ and $\psi'$: 
\begin{enumerate}
\item $\phi', \psi', 
\in \mathcal{L}_{*}$;

\item the branches of $\phi'$ and $\psi'$ starting at $x$ going up to the root consist only of Skeleton nodes.\footnote{The purpose of this restriction is to enforce preservation of non-empty joins by the term function $\phi'^{\mathbb{C}}$. The soundness of the rule is founded upon this and approximation of the argument $\gamma$ as the join of all closed elements below it. In the non-constructive setting of \cite{ConPal13} the same strategy is followed, except that the approximation is done by means of completely join-irreducibles. Since this can give rise to empty sets of approximants and hence empty joins, $+\vee$ is excluded in the analogous approximation rule in \cite{ConPal13}, as the join does not preserve empty joins coordinate-wise. In the present setting, the set of closed approximants is never empty, and hence this restriction may be dropped. Similar considerations apply to $-\wedge$.} 
    
\item for every node of the form $+\mu X.\psi(\overline{x}, X)$ or of the form $-\nu X.\phi(\overline{x}, X)$ in such branches, which is not in the scope of another binder,   all propositional variables and free fixed point variables in $\psi(\overline{x}, X)$ and $\phi(\overline{x}, X)$ must be among $\overline{x}$ and $X$; moreover,
    \begin{enumerate}
\item 
 the associated term function $\psi(\overline{x},X): (\mathbb{C}^{\delta})^\tau \times \mathbb{C}^{\delta} \rightarrow \mathbb{C}^{\delta}$ preserves $(\mathbb{C}^\tau \times \mathbb{C})$-targeted joins for all $\mu$-algebras
 $\mathbb{C}$ of the second kind; moreover  $\psi(\overline{x}, X)$ is required to be positive (negative) in $x_i$ if $\tau_i = 1$ ($\tau_i = \partial$), i.e.\ $\psi(\overline{x}, X)$ must be $\tau$-positive in $\overline{x}$;
\item 
the associated term function $\varphi(\overline{x}, X): (\mathbb{C}^{\delta})^\tau \times \mathbb{C}^{\delta}\rightarrow \mathbb{C}^{\delta}$ preserves $(\mathbb{C}^\tau \times \mathbb{C})$-targeted meets for all $\mathcal{L}_{\mathrm{LE}}$-algebras 
$\mathbb{C}$ of the second kind; moreover  $\phi(\overline{x}, X)$ is required to be positive (negative) in $x_i$ if $\tau_i = 1$ ($\tau_i = \partial$), i.e.\ $\phi(\overline{x}, X)$ must be $\tau$-positive in $\overline{x}$.
\end{enumerate}
\end{enumerate}

\begin{remark}
\begin{enumerate}
\item The approximation rules above, as stated, are sound both under admissible and under arbitrary assignments. However, their liberal application gives rise to topological complications in the proof of canonicity. Therefore, we will restrict the applications of approximation rules to nodes $!x$ giving rise to {\em maximal} skeleton branches.
    Such applications will be called {\em pivotal}. Also, executions of ALBA in which approximation rules are applied only pivotally will be referred to as {\em pivotal}.\label{pivotal:approx:rule:application}
\item In \cite{CoCr14}, approximation rules were formulated specifically for formulas having a fixed point binder as main connective. These rules had a substantially more cumbersome formulation than the one given above, which, modulo the restrictions about the preservation of targeted joins and meets, follows verbatim the approximation rules of \cite{ConPal13}. Moreover, the approximation rules \cite{CoCr14} could give rise to the splitting of the quasi-inequality into a set of quasi-inequalities, which is not the case of the present setting. This is thanks to  the fact that nominals and co-nominals are not interpreted as  completely join-primes (resp.\ meet-primes), but as closed and open elements, and this notion is compatible with products (cf.\ Section \ref{Sec:Prelim}, page \pageref{Comment:Compatibility}).

\end{enumerate}
\end{remark}

\paragraph{Residuation rules.} These rules operate on the inequalities in $S$, by rewriting a chosen inequality in $S$ into another inequality. For every $f\in \mathcal{F}$ and $g\in \mathcal{G}$, and any $1\leq i\leq n_f$ and $1\leq j\leq n_g$,

\begin{prooftree}
\AxiomC{$f(\phi_1,\ldots,\phi_i,\ldots,\phi_{n_f}) \leq \psi $}
\RightLabel{$\epsilon_f(i) = 1$}
\UnaryInfC{$\phi_i\leq f^\sharp_i(\phi_1,\ldots,\psi,\ldots,\phi_{n_f})$}
\AxiomC{$f(\phi_1,\ldots,\phi_i,\ldots,\phi_{n_f}) \leq \psi $}
\RightLabel{$\epsilon_f(i) = \partial$}
\UnaryInfC{$f^\sharp_i(\phi_1,\ldots,\psi,\ldots,\phi_{n_f})\leq \phi_i$}
\noLine\BinaryInfC{}
\end{prooftree}

\begin{prooftree}
\AxiomC{$\psi\leq g(\phi_1,\ldots,\phi_i,\ldots,\phi_{n_g})$}
\RightLabel{$\epsilon_g(i) = 1$}
\UnaryInfC{$g^\flat_i(\phi_1,\ldots,\psi,\ldots,\phi_{n_g})\leq \phi_i$}
\AxiomC{$\psi\leq g(\phi_1,\ldots,\phi_i,\ldots,\phi_{n_g})$}
\RightLabel{$\epsilon_g(i) = \partial$}
\UnaryInfC{$\phi_i\leq g^\flat_i(\phi_1,\ldots,\psi,\ldots,\phi_{n_g})$}
\noLine\BinaryInfC{}
\end{prooftree}
In a given system, each of these rules replaces an instance of the upper inequality with the corresponding instances of the two lower inequalities.

\paragraph{Ackermann rules}
The Ackermann rules are used for the crucial task of eliminating propositional variables from
quasi-inequalities.

\begin{prooftree}
 \AxiomC{$\exists p[\bigamp_{i = 1}^{n} \alpha_i \leq p \: \amp \: \bigamp_{j = 1}^{m} \beta_j(p) \leq \gamma_j(p)]$}
 \RightLabel{(RA)}
 \UnaryInfC{$\bigamp_{j = 1}^{m} \beta_j(\bigvee_{i = 1}^{n} \alpha_i / p) \leq \gamma_j(\bigvee_{i = 1}^{n} \alpha_i / p) $}
\end{prooftree}
subject to the restrictions that the $\alpha_i$ are $p$-free and syntactically closed, the $\beta_j$ are positive in $p$ and syntactically closed, while the $\gamma_j$ are negative in $p$ and syntactically open (cf.\ subsection below).

\begin{prooftree}
 \AxiomC{$\exists p[\bigamp_{i = 1}^{n} p \leq \alpha_i \: \amp \: \bigamp_{j = 1}^{m} \gamma_j(p) \leq \beta_j(p) ]$}
 \RightLabel{(LA)}
 \UnaryInfC{$\bigamp_{j = 1}^{m} \gamma_j(\bigwedge_{i = 1}^{n} \alpha_i / p) \leq \beta_j(\bigwedge_{i = 1}^{n} \alpha_i / p) $}
\end{prooftree}
subject to the restrictions that the $\alpha_i$ are $p$-free and syntactically open, the $\beta_j$ are positive in $p$ and syntactically open, while the $\gamma_j$ are negative in $p$ and syntactically closed.


\subsubsection{Syntactically (almost) open and closed formulas}
As mentioned in the introduction, when formulas from the extended language
$\Tm^+$ are interpreted in $\A^\delta$, an assignment $V$ will have
$V(\NOM)\subseteq K(\A^\delta)$ and $V(\CNOM) \subseteq O(\A^\delta)$.
For any assignment $V$, the assignment $V'$ is a $p$-variant (or, a $\nomj$-variant or
$\cnomm$-variant) of $V$ if $V'$ agrees with $V$ on all elements of $\mathsf{PROP}\cup\mathsf{NOM}\cup\mathsf{CNOM}$
except possibly at $p$ (respectively, at $\nomj$ or $\cnomm$). If so, we write $V' \sim_p V$
(respectively, $V' \sim_{\nomj} V$ or $V' \sim_{\cnomm} V$).
From this point on, we will use $V$ to denote both the assignment
$V : \mathsf{PROP} \cup \NOM \cup\CNOM \to \A^\delta$ and its unique homomorphic extension
$V: \Tm^+ \to \A^\delta$.

In the following definition we will use $f\in \mathcal{F}$ and $g\in \mathcal{G}$ to denote connectives of the original signature, and $h\in \mathcal{F}^+\setminus \mathcal{F}$ and $k\in \mathcal{G}^+\setminus \mathcal{G}$  to denote connectives of the expanded (`tense') language. To simplify notation, we will disregard the actual order of the coordinates, but keep track of their polarity. So, for instance we will write $f(\overline{\psi}, \overline{\phi})$ and $k(\overline{\phi}, \overline{\psi})$, where in both cases the coordinates are divided in two possibly empty arrays, the first (resp.\ second) of which contains the positive (resp.\ negative) coordinates.

\begin{definition}\label{def:synopensynclosed}
The \emph{syntactically open formulas} $\phi$ and \emph{syntactically closed formulas} $\psi$ are defined by simultaneous mutual recursion as follows:
\begin{eqnarray*}
\phi ::= \bot \mid \top \mid p \mid \cnomm \mid \phi_1 \wedge \phi_2 \mid \phi_1 \vee \phi_2 \mid g(\overline{\phi},\overline{\psi}) \mid f(\overline{\phi},\overline{\psi}) \mid k(\overline{\phi},\overline{\psi})
\mid \nu^{*} X. \phi\\
\psi ::= \bot \mid \top \mid p \mid \nomi \mid \psi_1 \wedge \psi_2 \mid \psi_1 \vee \psi_2 \mid f(\overline{\psi},\overline{\phi}) \mid g(\overline{\psi},\overline{\phi}) \mid h(\overline{\psi},\overline{\phi})
\mid \mu^{*} X. \psi
\end{eqnarray*}
where $p \in \mathsf{PROP}$, $\nomi \in \mathsf{NOM}$, and $\cnomm \in \mathsf{CNOM}$.

The \emph{syntactically almost open formulas} and \emph{syntactically almost closed formulas} are defined by adding $\mu^{*} X. \phi$ (resp.\ $\nu^{*} X. \psi$) to the recursive definition of $\phi$ (resp.\ $\psi$).

Informally, an $\Tm^+_*$-term is \emph{syntactically almost open} (resp.\ \emph{syntactically almost closed}) if, in it,
all occurrences of nominals and $h\in \mathcal{F}^+\setminus\mathcal{F}$ are negative (resp.\ positive), while all occurrences of co-nominals and $k\in \mathcal{G}^+\setminus\mathcal{G}$  are positive (resp.\ negative). If, in addition, all occurrences of $\mu^*$ are negative (resp.\ positive) and all occurrences of $\nu^*$ positive (resp.\ negative), the term is \emph{syntactically open} (resp.\ \emph{syntactically closed}).
\end{definition}

\subsection{Stage 3: Success, failure and output}

If stage 2 succeeded in eliminating all propositional variables from each system, the algorithm returns the conjunction of these purified quasi-inequalities, denoted by $\mathrm{ALBA}(\phi \leq \psi)$. Otherwise, the algorithm reports failure and terminates.

A \emph{tame run} of $\mu^*$-ALBA is one during which the approximation rules are applied only to formulas $\phi'(\gamma/!x)$ and $\psi'(\gamma/!x)$ such that no fixed point binder occurs in the branch from $x$ to the root of $\phi'$ and $\psi'$. 
By contrast, a \emph{proper run} of $\mu^*$-ALBA is one during which \emph{all} occurrences of fixed point binders lie along some branch ending with a subterm $\gamma$ which the application of an approximation rule extracts.
We say that a run of the algorithm $\mu^*$-ALBA \emph{succeeds} if all propositional variables are eliminated from the input inequality, $\phi \leq \psi$, and  denote the resulting set of pure quasi-inequalities by $\mathrm{ALBA}(\phi^*\le \psi^*)$. An inequality on which some run of $\mu^*$-ALBA succeeds is called a \emph{$\mu^*$-ALBA inequality}.

\end{document}

\section{Soundness of the approximation rule}
\subsection{The left arm}
Before we prove the soundness of the Approximation Rules ($\mu^\tau$-A-R) and ($\nu^\tau$-A-R), we require some lemmas regarding
preservation properties of operations.

\begin{lemma}\label{lem:heart} Let\, $\mathbf{L}$ be a lattice and $\tau$ an order-type over $m$. Let $f : \mathbf{L}^{m+1} \to \mathbf{L}$ be a term function such
that $f^{\mathbf{L}^\delta} : (\mathbf{L}^\delta)^\tau \times \mathbf{L}^\delta \to \mathbf{L}^\delta$ preserves
$\mathbf{L}^{\tau \oplus 1}$-targeted joins. Let $g: (\mathbf{L}^\delta)^\tau \to \mathbf{L}^\delta$ be the function given by $g(x_1,\ldots,x_m) = f(x_1,\ldots,x_m,\bot)$. Then $g$ preserves
$\mathbf{L}^\tau$-targeted joins.
\end{lemma}
\begin{proof}
\texttt{see the mu paper}
\end{proof}

We note that this result would also hold if we replaced $\bot$ with any $a \in L$, so long as $S$ is non-empty. To accommodate the case
that $S=\emptyset$, we must have $a = \bot$.

\begin{lemma}\label{lem:2.1equiv} Let $\A$ be a mu-algebra of the second kind and let $\tau$ be an order-type over $n$.
Let $\tau'$ be the order-type over $n+1$ defined by $\tau'=\tau \oplus 1$.
Suppose $f^{\A^\delta} : (\A^\delta)^{\tau'} \to \A^\delta$ is an $\mathcal{L}_{*}$ term function which
%
preserves $\A^{\tau'}$-targeted joins and is $\tau'$ positive.
Let $S \subseteq (\A^\delta)^{\tau}$ such that $\bigvee S \in \A^{\tau}$, i.e. $\bigvee S$ is an $\A^{\tau}$-targeted join. Then
$\LFPstar x. f^{\A^\delta}(\TJ S,x) = \bigvee \{\, \LFPstar x.f^{\A^\delta}(s,x)\mid s \in S\,\}$.
\end{lemma}

\begin{proof}
\texttt{see the mu paper}
\end{proof}

Before demonstrating that the Approximation Rule is sound we should point out that this rule is only ever
applied to a \emph{quasi-inequality} that is the result of an application of the First Approximation Rule (FA).
That is, ($\mu^\tau$-A-R) and ($\nu^\tau$-A-R) are each applied to an inequality which forms \emph{part of}
the antecedent of an implication.
When demonstrating the soundness of the rule ($\mu^\tau$-A-R) it is therefore sufficient
to show that the inequality above the line and the inequality below the line are valid under
assignments
which agree everywhere except at some nominal which does not occur in the consequent of the quasi-inequality.

Some remarks are in order:
Firstly, the reason that we need to use $\mu^*$ and $\nu^*$ in the rules ($\mu^\tau$-A-R) and ($\nu^\tau$-A-R), respectively, is shown by
the proof of Lemma~\ref{lem:2.1equiv}. If we were to use $\mu x.f^{\A^\delta}(\TJ S,x)=\bigwedge \{\,a \in A^\delta \mid f^{\A^\delta}(\bigvee S,a) \le a\,\}$
in line (1),
we would then have only $(1) \le (2)$ (as $A \subseteq A^\delta$). Thus if we formulated ($\mu^\tau$-A-R) and Lemma~\ref{lem:2.1equiv} with
$\mu$ instead of $\mu^*$, we would not have the equality in Lemma~\ref{lem:2.1equiv} and thus would not be able to show the invariance of admissible validity under ($\mu^\tau$-A-R).

Secondly, note that we cannot generalize Proposition \ref{prop:approxsound} to mu-algebras of the first kind. This is because its proof depends on Lemma \ref{lem:heart} which crucially makes use of ordinal-unfolding of fixed-points, which is only available in mu-algebras of the second kind.

\texttt{Mention that inner formulas are the ones where the targeted join preservation is guaranteed.}

In the formulation of the approximation rules ($\mu^{\tau}$-A-R) and ($\nu^{\tau}$-A-R) we required the term functions $\psi(\overline{x},X)$ and $\phi(\overline{x},X)$
to preserve, respectively, targeted $\A^\tau$-joins and targeted $\A^\tau$-meets.
The inner formulas (introduced in \cite[Section 4]{muALBA}) are a syntactically specified classes of formulas, the term functions of which satisfy these properties. Being an inner formula is thus an effectively checkable sufficient condition for the applicability of the rules ($\mu^{\tau}$-A-R) and ($\nu^{\tau}$-A-R).\footnote{In a related line of research, Fontaine \cite{Fontaine08} gives a syntactic characterization of the \emph{continuous} formulas of the (classical) modal mu-calculus. Continuity is shown to be equivalent to Scott-continuity, which is easily seen to imply the preservation of non-empty joins. It does not, however, imply the preservation of the empty join, as can be seen by considering the formula $\top \vee p$, for example.}

\begin{definition}\label{IF:Definition}
Let $\oy, \oz \subseteq \mathsf{PHVAR}$ and $\oX \subseteq \mathsf{FVAR}$ be tuples of variables which are pairwise different in the union of the their underlying sets. Let $\tau$ be an order-type on $\ox = \oy \oplus \oX$. The {\em $\tau$-$\Box$ and $\tau$-$\Diamond$ $(\ox, \oz)$-inner formulas} (\emph{$(\ox, \oz)$-\ifBox{\tau}} and \emph{$(\ox, \oz)$-\ifDia{\tau}}), the free variables of which are contained in $(\ox, \oz)$, are given by the following simultaneous recursion (for the sake of readability, the parameters $\ox$ and $\oz$ are omitted):

\begin{center}
\begin{tabular}{lcccccccccccccc}
\ifBox{\tau}$\ni \phi$ &$\!\!::=\!\!$ &$x_i$ &$\!\!\mid\!\!$ &$\Box \phi$ &$\!\!|\!\!$ &$\phi_1 \wedge \phi_2$ &$\!\!|\!\!$ &$\nu^* Y. \phi' $ &$\!\!|\!\!$ &$\pi \rightarrow \phi$ &$\!\!|\!\!$ &$\pi \vee \phi$ &$\!\!|\!\!$ &$\psi^c \rightarrow \pi$\\
\ifDia{\tau}$\ni \psi$ &$\!\!::=\!\!$ &$x_i$ &$\!\!\mid\!\!$ &$\Diamond \psi$ &$\!\!|\!\!$ &$\psi_1 \vee \psi_2$ &$\!\!|\!\!$ &$\mu^* Y. \psi' $ &$\!\!|\!\!$ &$\psi - \pi$ &$\!\!|\!\!$ &$\pi \wedge \psi$ &$\!\!|\!\!$ &$\pi - \phi^c$\\
\end{tabular}
\end{center}\marginnote{Of course, adapt this into the LE setting.}
\noindent where
\begin{enumerate}
\item $\tau_i = 1$ in the base of the recursion,\label{IF:Base}
\item $\pi$ is $\pi(\oz)\in \mathcal{L}_*$ (specifically, $\pi(\oz)$ contains none of the variables in $\ox$ or $\oX$),\label{IF:Pi}
\item $\phi' = \phi'(\oy \oplus \oX',\oz)$ and $\psi' = \psi'(\oy \oplus \oX',\oz)$ are \ifBox{\tau'} and \ifDia{\tau'}, respectively, with $\oX' = \oX \oplus Y$ and $\tau' = \tau \oplus 1$,\label{IF:FixPoint}
\item $\psi^c\in (\ox, \oz)$-\ifDia{\tau^{\partial}} and $\phi^c\in (\ox, \oz)$-\ifBox{\tau^{\partial}}.\label{IF:Reverse}
\item All other formulas have their free variables among $(\ox, \oz)$.
\end{enumerate}
\end{definition}

\noindent The key fact about $(\ox, \oz)$-\ifBox{\tau} and $(\ox, \oz)$-\ifDia{\tau} formulas is the following:

\begin{lemma}\label{IF:Are:Adjoint}
Let $\A$ be a mu-algebra of the second kind.

\begin{enumerate}
\item The term function associated with any
\ifDia{\tau} formula $\psi (\ox, \oz)$
(resp., \ifBox{\tau} formula $\phi (\ox, \oz)$ )
will preserve $\A^\tau$-targeted joins (resp., targeted $\A^\tau$ meets) as a map from $(\A^{\delta})^\tau$ to $\A^{\delta}$
%
when the variables $\oz$ are fixed.


\item The term function associated with any
\ifDia{\tau} formula $\psi (\ox, \oz)$
(resp., \ifBox{\tau} formula $\phi (\ox, \oz)$ ) that does not contain any occurrences of either $\mu^*$ or $\nu^*$
will be completely join-preserving (resp., completely meet-preserving) as a map from $(\A^{\delta})^\tau$ to $\A^{\delta}$ when the variables $\oz$ are fixed.
\item If the $\oga$ are constant $\mathcal{L}_*^+$ sentences, then the term function associated with $\psi (\ox, \oga/ \oz)$ (resp., $\phi (\ox, \oga/ \oz)$)
will preserve $\A^\tau$-targeted joins (resp., targeted $\A^\tau$ meets) as a map from $(\A^{\delta})^\tau$ to $\A^{\delta}$.
\item Every \ifDia{\tau} formula $\psi (\ox, \oz)$ and every \ifBox{\tau} formula $\phi (\ox, \oz)$ is $\tau$-positive in $\ox$.
\end{enumerate}
\end{lemma}
\begin{proof}
Let $\tau$ be an order-type over $n$ and let $S \subseteq (A^\delta)^n$ be such that $\TJ S$ is a $\A^\tau$-targeted join. That is, $\TJ^{(\A^\delta)^\tau} S \in \A^\tau$.
We note that this also implies that $\TM S$ is a $\A^{\tau^\partial}$-targeted meet since $\TJ^{(\A^\delta)^\tau} S = \TM^{(\A^{\delta^\partial})^\tau} S$.

Claim (3) is an immediate consequence of claim (1) and claim (4) follows by a straightforward simultaneous induction on $\psi$ and $\phi$. We proceed by simultaneous induction on $\psi$ and $\phi$ to prove claims (1) and (2). The base cases are clear, as are the cases for $\Box \phi$, $\Diamond \psi$, $\phi_1\wedge \phi_2$, and $\psi_1 \vee \psi_2$.
If $\psi$ is of the form
$\mu^*Y.\psi'$ with $\psi(\oy \oplus \oX,\oz)$ and $\psi'(\oy\oplus\oX\oplus Y,\oz)$ then the result follows from the induction hypothesis and the result of Lemma~\ref{lem:2.1equiv}. That is,
$$\psi(\TJ S,\oz) = \mu^*Y.\psi'(\TJ S, Y,\oz)=\bigvee \{\, \mu^*Y. \psi'(s,Y,\oz) \mid s \in S \,\}=\bigvee \{\,\psi(s,\oz)\mid s \in S\,\}.$$
The order dual of Lemma~\ref{lem:2.1equiv} can be used for the case of $\phi(\oy\oplus\oX,\oz)$ being of the form $\nu^*Y. \phi'(\oy\oplus\oX\oplus Y,\oz)$.

The final case that we will consider is if $\psi$ is of the form $\pi - \phi^c$ where $\phi^c \in (\ox,\oz)$--\ifBox{\tau^\partial}. (The other remaining cases follow similarly.)
Then
\begin{align*}
\psi(\TJ^\tau S,\oz) &= \pi(\oz) - \phi^c(\TJ^\tau S)
 = \pi(\oz) - \phi^c(\TM^{\tau^\partial}S,\oz) \\
& = \pi(\oz) - \bigwedge \{\, \phi^c(s,\oz) \mid s \in S \,\}\\
& = \bigvee \{\, \pi(\oz)- \phi^c(s,\oz) \mid s \in S \,\} = \bigvee\{\,\psi(s,\oz)\mid s \in S\,\}.
\end{align*}
\end{proof}

For our purposes then, the most important consequence of this lemma is the fact that ($\mu^{\tau}$-A-R) and ($\nu^{\tau}$-A-R) are respectively applicable to formulas of the form $\psi (\ophi/\ox, \oga/ \oz)$ and $\phi (\opsi / \ox, \oga/ \oz)$, where $\psi (\ox, \oz)$ is \ifDia{\tau}, $\psi (\ox, \oz)$ is \ifBox{\tau} and the $\oga$ are constant $\mathcal{L}_*^+$ sentences. In particular, ($\mu^{\tau}$-A-R) and ($\nu^{\tau}$-A-R) are applicable to formulas $\psi (\ophi/\ox)$ and $\phi (\opsi / \ox)$, where $\psi (\ox, \oz)$ is \ifDia{\tau} and $\phi (\ox, \oz)$ is \ifBox{\tau} with $\oz$ the empty tuple.

\texttt{skeletons are inner formulas.}

\texttt{Here we put the distribution lemma and the soundness proof.}

\subsection{The right arm}

The next two lemmas and the proposition which follows are needed to convert formulas of the form $\mu^*X.\phi(X)$ back to $\mu X.\phi(X)$.
The first of the lemmas is proved in~\cite{muALBA}.

\begin{lemma}\label{lem:2.1orig}{\upshape \cite[Lemma 2.1]{muALBA}} Let $\mathbf{L}$ and $\mathbf{M}$ be complete lattices and
$G : M \times L \to L$. Let $\mu y.G(-,y) : \mathbf{M} \to \mathbf{L}$ be given for $a \in M$ by
$a \mapsto \bigwedge \{\, x \in L \mid G(a,x) \le x\,\}$.

If $G$ is completely join-preserving, then $\mu y.G(-,y) : \mathbf{M} \to \mathbf{L}$ is defined everywhere on $M$ and is
completely join-preserving.
\end{lemma}

\begin{lemma}\label{lem:droppingstar} Let $\A$ be a mu-algebra of the second kind and let $\tau$ be an order-type over $n$. Let
$\psi(\overline{x},X) \in \Tm_1$ such that $\psi(\overline{x},X)$ is completely
$\bigvee$-preserving in $(\overline{x},X) \in (\A^\delta)^\tau \times \A^\delta$.
If $\A^\delta \models [(\nomi \le \mu^* X. \psi(\overline{\nomj_i}^\tau,X)
\amp \nomj^{\tau_i} \le^{\tau_i} \overline{\phi}_i)\Rightarrow \nomi \le \cnomm]$
for all $i=1,\ldots,n$,
then
$\A^\delta \models [\nomi \le \mu X. \psi(\overline{\phi},X) \Rightarrow \nomi \le \cnomm]$.
\end{lemma}\marginnote{For this lemma, we need to adapt into the non-distributive setting, since we do not have the double approximation anymore.}
\begin{proof}
Assume that $\A^\delta \models [(\nomi \le \mu^* X. \psi(\overline{\nomj_i}^\tau,X) \amp \nomj^{\tau_i} \le^{\tau_i} \overline{\phi}_i)\Rightarrow \nomi \le \cnomm]$ and suppose $V(\nomi) \le V(\mu X.\psi(\overline{\phi},X))$. 
By the join-density of $\jty(\A^\delta)$ we have that
$\overline{\varphi}=\bigvee\{\, \overline{j} \mid \overline{j} \in \jty((\A^\delta)^\tau), \overline{j} \le \overline{\phi}\,\}$. Thus
$\mu X.\psi(\overline{\phi},X)=\mu X. \psi( \bigvee\{\,\overline{j}\mid \overline{j}\le \overline{\phi}\,\},X)$.
Now by Lemma~\ref{lem:2.1orig} with $M=(\A^\delta)^\tau$ and $L=\A^\delta$ we get
$\mu X.\psi(\overline{\phi},X)=\bigvee\{\, \mu X.\psi(\overline{j},X) \mid \overline{j} \le \overline{\phi}\,\}$.
Since $V(\nomi) \in \jty(\A^\delta)$, it is completely join-prime and so
$V(\nomi) \le \mu X.\psi(\overline{\phi},X)=\bigvee\{\, \mu X.\psi(\overline{j},X) \mid \overline{j} \le \overline{\phi}\,\}$
implies that there exists some $\overline{j_0}$ with $\overline{j_0} \le \overline{\phi}$ such that
$V(\nomi) \le \mu X.\psi(\overline{j_0},X) =\bigwedge \{\,a \in A^\delta \mid \psi(\overline{j_0},a)\le a \,\}
\le \bigwedge \{\, a \in A \mid \psi(\overline{j_0},a)\,\}= \mu^* X.\psi(\overline{j_0},X)$.
Thus by the assumption, $V(\nomi) \le V(\cnomm)$.
\end{proof}

\section{Tame canonicity and proper canonicity}

\section{Success}

For tame canonicity, just take the mu paper, and adapt the approximation steps, which should be fine.

For proper canonicity,

\bibliographystyle{aiml16}
\bibliography{aiml16}

\end{document}*

%% file: introduction.tex
The present contribution lies at the crossroads of at least three active lines of research in nonclassical logics: the one investigating the semantic and proof-theoretic environment of fixed point expansions of logics algebraically captured by varieties of (distributive) lattice expansions \cite{Baelde,Hartonas98,santocanale05,BH12,Gavazzo}; the one investigating constructive canonicity for intuitionistic and substructural logics \cite{GhMe97,Suzuki11a}; the one uniformly extending  the state-of-the-art in Sahlqvist theory to families of nonclassical logics, and applying it to issues both semantic and proof-theoretic \cite{CoGhPa14,ConPalSou12,ConPal12,CoFoPaSo15,CoRo14,PaSoZh15b,PaSoZh15a,CoPaSoZh15,CoPaZh15a,ConPal13,GrMaPaTzZh15,MaZh16,FrPaSa14}, known as `unified correspondence'. 

We prove the algorithmic canonicity of two classes of $\mu$-inequalities in a constructive meta-theory of normal lattice expansions. This result simultaneously generalizes Conradie and Craig's canonicity results for $\mu$-inequalities based on a bi-intuitionistic bi-modal language \cite{CoCr14}, and Conradie and Palmigiano's constructive canonicity for inductive inequalities \cite{CoPaZh16b} (restricted to normal lattice expansions to keep the page limit). Besides the greater generality, the unification of these strands smoothes the existing proofs for the canonicity of $\mu$-formulas and inequalities. Specifically, the two  canonicity results proven in \cite{CoCr14}, (namely, the tame and proper canonicity, cf.\ Section \ref{Sec:TwoKindsCanon}) fully generalize to the constructive setting and  normal LEs. Remarkably, the rules of the algorithm ALBA used for this result have exactly the same formulation as those of \cite{CoPaZh16b}, with no additional rule added specifically to handle the fixed point binders. Rather, fixed points are accounted for by certain restrictions on the application of the rules, concerning the order-theoretic properties of the term functions associated with the formulas to which the rules are applied.   

Applications of these results include the formalization of common knowledge-type processes of social interaction. For instance, in ongoing work \cite{categorization} we are exploring the formalization of processes giving rise to categorization systems  (where categories arise, as in Formal Concept Analysis, as Galois-stable sets from a given polarity) agreed upon by a whole community.

\paragraph{Structure of the paper.}
In Section \ref{Sec:Prelim}, we collect the needed preliminaries.  In Section \ref{sec:synclass}, we define the classes of mu-inequalities to which the canonicity results apply. In Section \ref{Sec:TwoKindsCanon}, we introduce and expand on the two notions of canonicity mentioned above.  In Section \ref{sec:constr:mu}, we outline the constructive and general lattice version of the algorithm ALBA which is the main tool for the two canonicity results (more details are given in the appendix).  In Section \ref{Sec:Main:Results}, we state our main results. Due to space constraints we do not include proofs, but these may be found online in the an expanded version of the present paper \cite{long-version}.

%% file: sketch_alessandra.bbl
\begin{thebibliography}{10}
\expandafter\ifx\csname url\endcsname\relax
  \def\url#1{\texttt{#1}}\fi
\expandafter\ifx\csname urlprefix\endcsname\relax\def\urlprefix{URL }\fi
\newcommand{\enquote}[1]{``#1''}

\bibitem{AmKwMe95}
Ambler, S., M.~Kwiatkowska and N.~Measor, \emph{Duality and the completeness of
  the modal $\mu$-calculus}, Theoretical Computer Science \textbf{151} (1995),
  pp.~3--27.

\bibitem{Baelde}
Baelde, D., \emph{Least and greatest fixed points in linear logic}, ACM
  Transactions on Computational Logic (TOCL) \textbf{13} (2012), p.~2.

\bibitem{BH12}
Bezhanishvili, N. and I.~Hodkinson, \emph{Sahlqvist theorem for modal fixed
  point logic}, Theoretical Computer Science \textbf{424} (2012), pp.~1--19.

\bibitem{CoCr14}
Conradie, W. and A.~Craig, \emph{Canonicity results for mu-calculi: An
  algorithmic approach}, Journal of Logic and Computation  (forthcoming).
\newline\urlprefix\url{http://arxiv.org/abs/1408.6367}

\bibitem{long-version}
Conradie, W., A.~Craig, A.~Palmigiano and Z.~Zhao, \emph{Constructive
  canonicity for lattice-based fixed point logics}, ArXiv preprint .

\bibitem{CoFoPaSo15}
Conradie, W., Y.~Fomatati, A.~Palmigiano and S.~Sourabh, \emph{Algorithmic
  correspondence for intuitionistic modal mu-calculus}, Theoretical Computer
  Science \textbf{564} (2015), pp.~30--62.

\bibitem{categorization}
Conradie, W., S.~Frittella, A.~Palmigiano, M.~Piazzai, A.~Tzimoulis and
  N.~Wijnberg, \emph{Reasoning about categories with generalized {K}ripke
  frames}, in preparation .

\bibitem{CoGhPa14}
Conradie, W., S.~Ghilardi and A.~Palmigiano, \emph{Unified correspondence}, in:
  A.~Baltag and S.~Smets, editors, \emph{Johan van Benthem on Logic and
  Information Dynamics},  Outstanding Contributions to Logic  \textbf{5},
  Springer International Publishing, 2014 pp. 933--975.

\bibitem{CoPaZh16b}
Conradie, W. and A.~Palmigiano, \emph{Constructive canonicity of inductive
  inequalities} Submitted.

\bibitem{ConPal12}
Conradie, W. and A.~Palmigiano, \emph{Algorithmic correspondence and canonicity
  for distributive modal logic}, Annals of Pure and Applied Logic \textbf{163}
  (2012), pp.~338 -- 376.

\bibitem{ConPal13}
Conradie, W. and A.~Palmigiano, \emph{Algorithmic correspondence and canonicity
  for non-distributive logics}, Journal of Logic and Computation
  (forthcoming).

\bibitem{ConPalSou12}
Conradie, W., A.~Palmigiano and S.~Sourabh, \emph{Algorithmic modal
  correspondence: {S}ahlqvist and beyond} Submitted.

\bibitem{CoPaSoZh15}
Conradie, W., A.~Palmigiano, S.~Sourabh and Z.~Zhao, \emph{Canonicity and
  relativized canonicity via pseudo-correspondence: an application of {ALBA}}
  Submitted.

\bibitem{CoPaZh15a}
Conradie, W., A.~Palmigiano and Z.~Zhao, \emph{Sahlqvist via translation}
  Submitted.

\bibitem{CoRo14}
Conradie, W. and C.~Robinson, \emph{On {S}ahlqvist theory for hybrid logic},
  Journal of Logic and Computation  (2015).

\bibitem{DunnGP05}
Dunn, J.~M., M.~Gehrke and A.~Palmigiano, \emph{Canonical extensions and
  relational completeness of some substructural logics}, Journal Symbolic Logic
  \textbf{70} (2005), pp.~713--740.

\bibitem{FrPaSa14}
Frittella, S., A.~Palmigiano and L.~Santocanale, \emph{Dual characterizations
  for finite lattices via correspondence theory for monotone modal logic},
  Journal of Logic and Computation  (forthcoming).

\bibitem{Gavazzo}
Gavazzo, F., \emph{Investigations into linear logic with fixed-point operators}
   (2015), iLLC MoL Thesis.

\bibitem{GeHa01}
Gehrke, M. and J.~Harding, \emph{Bounded lattice expansions}, Journal of
  Algebra \textbf{238} (2001), pp.~345--371.

\bibitem{GhMe97}
Ghilardi, S. and G.~Meloni, \emph{Constructive canonicity in non-classical
  logics}, Annals of Pure and Applied Logic \textbf{86} (1997), pp.~1--32.

\bibitem{GrMaPaTzZh15}
Greco, G., M.~Ma, A.~Palmigiano, A.~Tzimoulis and Z.~Zhao, \emph{Unified
  correspondence as a proof-theoretic tool}, Journal of Logic and Computation
  (forthcoming).

\bibitem{Hartonas98}
Hartonas, C., \emph{Duality for modal $\mu$-logics}, Theoretical Computer
  Science \textbf{202} (1998), pp.~193--222.

\bibitem{MaZh16}
Ma, M. and Z.~Zhao, \emph{Unified correspondence and proof theory for strict
  implication}, Journal of Logic and Computation  (forthcoming).

\bibitem{PaSoZh15a}
Palmigiano, A., S.~Sourabh and Z.~Zhao, \emph{J\'onsson-style canonicity for
  {ALBA}-inequalities}, Journal of Logic and Computation  (2015).

\bibitem{PaSoZh15b}
Palmigiano, A., S.~Sourabh and Z.~Zhao, \emph{Sahlqvist theory for impossible
  worlds}, Journal of Logic and Computation  (forthcoming).

\bibitem{santocanale05}
Santocanale, L., \emph{Completions of $\mu$-algebras}, in: \emph{Logic in
  Computer Science, 2005. LICS 2005. Proceedings. 20th Annual IEEE Symposium
  on}, IEEE, 2005, pp. 219--228.

\bibitem{Suzuki11a}
Suzuki, T., \emph{Canonicity results of substructural and lattice-based
  logics}, The Review of Symbolic Logic \textbf{4} (2011), pp.~1--42.

\end{thebibliography}
